\documentclass[12pt,a4paper]{article}

\usepackage[english]{babel}
\usepackage{indentfirst}
\usepackage{misccorr}
\usepackage{graphicx}
\usepackage{amsmath}
\usepackage{amssymb}
\usepackage{amsthm}
\usepackage{hyperref}
\usepackage{color}

\newif\ifhideru
\hiderufalse

\newif\ifhideen
\hideentrue 

\newlength{\wdth}

\newtheorem{theorem}{Theorem}

\newtheorem{rTheorem}{Теорема}
\newtheorem{rLemma}{Лемма}
\newtheorem{rProposition}{Предложение}
\newtheorem{rDefinition}{Определение}
\newtheorem{rRemark}{Замечание}

\newtheorem{rEx}{Пример}

\newtheorem{Lemma}{Lemma}
\newtheorem{Proposition}{Proposition}
\newtheorem{Definition}{Definition}
\newtheorem{Remark}{Remark}

\newtheorem{Ex}{Example}

\ifhideen
\title{{\normalsize\tt\hfill\jobname.tex}\\
On nonlinear weak  law of large numbers
}

\author{A.T. Akhmiarova\footnote{M.V. Lomonosov Moscow State University, Leninskie Gory, Moscow, 119991, Russian Federation \& Institute for Information Transmission Problems, 19 Bolshoy Karetny per., build. 1, 127051, Moscow, Russian Federation; email: lili-r01@yandex.ru}, \; A.Yu. Veretennikov\footnote{Institute for Information Transmission Problems, 19 Bolshoy Karetny per., build. 1, 127051, Moscow, Russian Federation \& RUDN University, 6 Miklukho Maklaya St., 117198, Moscow, Russian Federation, email: ayv@iitp.ru}}

\fi

\ifhideru
\title{{\normalsize\tt\hfill\jobname.tex}\\
О нелинейном слабом законе больших чисел
}

\author{А.Т. Ахмярова\footnote{Московский государственный университет им. М.В. Ломоносова, Ленинские горы, Москва, 119991, Российская Федерация \& Институт проблем передачи информации; email: lili-r01@yandex.ru}, \; А.Ю. Веретенников\footnote{Институт проблем передачи информации, 19 Bolshoy Karetny per., build. 1, 127051, Moscow, Russian Federation \& Российский университет дружбы народов, 6 Miklukho Maklaya St., 117198, Moscow, Russian Federation email: alexander.veretennikov2011@ya.ru}}
\fi

\date{}

\begin{document}
\maketitle

\ifhideen
\begin{abstract}
\noindent
A new  version of a weak nonlinear law of large numbers proposed. The existence of the first moment for any summand is not assumed. The assumption of independence is understood in the nonlinear sense, and may be further a little relaxed. 
\medskip

\noindent
Keywords: Nonlinear weak law of large numbers; nonlinear independence

\medskip

\noindent
MSC2020: 60F15

\end{abstract}
\fi

\ifhideen
\section{Introduction}
The basis of this paper is Kolmogorov's weak  law of large numbers (LLN in what follows)  \cite{Kolmogorov0,Kolmogorov} for i.i.d.r.v's $X_1, X_2, \ldots$ under the standing condition (here the notation $ \bar F_{X_1}(t) := 1- F_{X_1}(t)$ is in use)
\begin{equation}\label{ANKcondition}
t(F_{X_1}(-t) + \bar F_{X_1}(t))\to 0, \quad t\to\infty.
\end{equation}
The goal of this paper is to enchance its extension \cite[Theorem 3]{Akhmiarova-V} for non-identically distributed weakly dependent r.v's to the setting of a sublinear expectation. The latter is understood as proposed by S. Peng \cite{Peng1, Peng2}, although, without any link to BSDEs or likewise topics, just by using the basic definition of sublinear expectations. It may also be noted that apparently the first version of a generalized expectation  which loses certain linear properties was discussed in \cite[Ch.5.2]{Kolmogorov}, where it was denoted by $M^\circ$.

Various versions of weak (and in some of them also strong) nonlinear LLN may be found in 
\cite{Peng1, Peng2}, 
\cite{LiSun}; see also the references therein.  
In \cite{Fang}, \cite{Song} there are also estimates of the rate of convergence. 

As it was already said, nonlinear strong LLN results also exist, see, for example, \cite{ChengHu2018} where a version of such a theorem under the ``controlled first moment condition'' has been established along with a review of some other nonlinear SLLN theorems. We do not touch this subject here, but refer to certain preliminary results about nonlinear expectation in it.

Let us mention one feature of this work. Instead of ``usual'' truncations of r.v's at some level, we apply some ``smoothed version'' of this truncation where instead of indicators certain continuous functions with compact support are used. In principle, this is not necessary. However, this allows references to already established representations of nonlinear expectations, while with indicators such representations should have been re-proved. The changes in the calculi are minimal.

The present paper consists of four sections: 1 -- this Introduction, 2 -- Basics of  nonlinear expectations, 3 -- the Setting and main result, 4 -- the Proof (of theorem \ref{thm1}). 
\fi

\ifhideen
\section{Basics of nonlinear expectations}
The main ref on \cite{Peng1}, \cite{Lebedev2}, 
and \cite{ChengHu2018}. Firstly, let us describe the space of random variables on which the operation of nonlinear expectation and nonlinear probability (capacity) will be defined.

\begin{Definition}\label{def1}
Let $(\Omega, \mathcal{F})$ be a measurable space and $\cal H$ be a subset of all random variables on it such that: 
\begin{itemize}
    \item $I_A \in \cal H$ for any $A \in \mathcal{F}$,
    \item if $X_1, X_2, \dots, X_n \in\cal H$, then $\varphi(X_1, X_2, \dots, X_n) \in \cal H$ for any $\varphi \in C_{\ell,Lip}(\mathbb{R}^n)$, 
where $C_{\ell,Lip}(\mathbb{R}^n)$ is the space of locally Lipschitz functions satisfying the inequality 
\begin{equation}
    |\varphi(x) - \varphi(y)| \leq C(1 + |x|^m + |y|^m)|x - y|, \quad \forall x, y \in \mathbb{R}^n,
\end{equation}
with some $C > 0$ and $m \in \mathbb{N}$; both constants here may depend on $\varphi$. 

\end{itemize}
\end{Definition}

\begin{Definition}[Upper sublinear expectation $\mathbb E$]
Sublinear expectation $\mathbb{E}$ on $H$ is a functional $\mathbb{E}: H \to \overline{\mathbb{R}}:=[-\infty, \infty]$ satisfying the following properties: for  any $X, Y \in H$ it holds, 
\begin{enumerate}
    \item Monotonicity: $X \geq Y$ implies $\mathbb{E}[X] \geq \mathbb{E}[Y]$.
    \item Preserving of constants: $\mathbb{E}[c] = c$, $\forall c \in \mathbb{R}$.
    \item Positive homogeneity: $\mathbb{E}[\lambda X] = \lambda \mathbb{E}[X]$, $\forall \lambda \geq 0$.
    \item Subadditivity: $\mathbb{E}[X+Y] \leq \mathbb{E}[X] + \mathbb{E}[Y]$.
\end{enumerate}
\end{Definition}

\begin{Definition}[Lower sublinear  expectation $\cal E$]
The conjugate to $\mathbb E$, ``lower'' expectation $\cal E$ on $H$ is defnied by the formula, 
$$
{\cal E} X \equiv {\cal E} (X) := - \mathbb E (-X).
$$
\end{Definition}

\begin{Remark}
Sigma-algebra ${\cal F}$ is assumed to be completed\footnote{This completion is not necessary for a weak LLN.} by all sets $A$  (called ``polar sets'') satisfying $\mathbb E 1(A) = 0$. 
\end{Remark}

\begin{Ex}[\cite{Peng1}]\label{Ex1}
Let ${\cal P}$ be a family of probability measures on $(\Omega, {\cal F})$. Then $\mathbb E X := \sup_{\mathsf P\in {\cal P}} \mathsf E^{\mathsf P} X$ is an (upper) sublinear expectation in the sense of the definition~\ref{def1}.
\end{Ex} 
This example is a case in point, as in various applications there might be more than one probability measure and we may not know which one should be used, e.g., in incomplete markets in financial mathematics.

\begin{Definition}[Upper and lower capacities]
The function $V$ on $\cal F$ generated by the operator $\mathbb E$ is called upper capacity; the function $v$ on $\cal F$ is called lower capacity:
$$
V(A) := \mathbb E 1(A), \quad v(A) := {\cal E} 1(A).
$$
\end{Definition}

\begin{Definition}\label{def-indep}
Pairwise independence for r.v. $X_1,X_2$: for any 
test functions $\varphi_1, \varphi_2 \in C_{\ell, Lip}$, 
$$
\mathbb E \varphi_1(X_1)\varphi_2(X_2) = \mathbb E \varphi_1(X_1) \mathbb E\varphi_2(X_2).
$$
\end{Definition}

\medskip

The important Markov's and Bienaimé -- Chebyshev's inequalities hold.  
\begin{Lemma}[Markov’s and Bienaimé -- Chebyshev's inequalities, \cite{Peng1}, Lemma 2.2 \cite{ChengHu2018}, Proposition 2.1 \cite{ChenWuLi}]\label{Markov} 
Let $f(x) > 0$ be a nondecreasing function on $\mathbb R$.
Then for any $x$,
$$
V (X \ge x) \le \frac{\mathbb E [f(X)]}{f(x)}, \quad 
v(X \ge x) \le \frac{{\cal E}[f(X)]}{f(x)}.
$$
Also, for any $c>0$ and any r.v. $X$, 
$$
V (|X - \mathbb E X|  \ge c) \le \frac{\mathbb E |X - \mathbb E X| ^2}{c^2}, \quad 
v (|X - {\cal E} X|  \ge c) \le \frac{{\cal E} |X - {\cal E} X| ^2}{c^2}.
$$
\end{Lemma}

Another important tool is the following result.
\begin{Proposition}[Lemma 1.3.4 \cite{Peng1}]\label{pro-sadd}
For any r.v. $X \in H$ there exists a family of probability measures (sigma-additive) ${\cal P}= (\mathsf P^\theta, \theta\in \Theta)$ with some parametric set $\Theta$, so that $\mathbb E \varphi(X) = \sup_{\theta\in \Theta} \mathsf E^\theta \varphi(X)$ for any $\varphi \in C_{\ell,Lip}(\mathbb{R}^n)$. 
\end{Proposition}

\begin{Remark}\label{rem-fadd}
Notice that for a collection -- for more than one -- of random variables $\cal H$, generally speaking, only the existence of a family of {\bf finitely additive} measures may be guaranteed with the analogous property $\mathbb E \varphi(X) = \sup_{\theta\in \Theta} \mathsf E^{\mathsf P^\theta} \varphi(X)$, see \cite[Theorem 1.2.4]{Peng1}. However, this is not fully sufficient for our method of proof in what follows.
\end{Remark}
\fi

\ifhideen
\section{The setting and main result} 
Let us consider the sequence of random variables $(X_1, X_2, \ldots)$ on the nonlinear probability quadruple $(\Omega, {\cal F}, {\cal H}, \mathbb E)$.

Let $\chi_n(x)$ be a smooth (Lipschitz) bounded function equal to $1$ for $|x|\le n$, and to zero for  $|x| \ge n+1$, and taking values from $[0,1]$ for $n\le |x| \le n+1$, and let $\tilde\chi_n(x) := 1-\chi_n(x)$. Notice that $\tilde\chi_n(x) \le 1(|x|>n)$. In this paper, instead of standard truncations we will use functions $\chi_n(x)$. Let 
$$
Y_{n,k} = X_k\chi_n(|X_k|), \quad 
\tilde \mu^+_{n,k} :=  \mathbb E Y_{n,k}, \quad  
\tilde \mu^-_{n,k} :=  {\cal E} Y_{n,k},
$$
and
$$
\bar \mu_n := 
\frac1{n} \sum_{k=1}^n \tilde\mu^+_{n,k}, \quad 
\underline \mu_n := 
\frac1{n} \sum_{k=1}^n \tilde\mu^-_{n,k}.
$$

Let 
$$
S_n = \sum_{k=1}^n X_k, \quad S'_n = \sum_{k=1}^n Y_{n,k}.
$$
By virtue of the proposition \ref{pro-sadd}, for each $k$ there exists a parametric set $\Theta_{k}$ and a family of (sigma-additive) probability measures $(\mathsf P^\theta, \theta\in \Theta_{k})$ such that $\mathbb E \varphi(X_{k}) = \sup_{\theta\in \Theta_{k}} \mathsf E^{\theta} \varphi(X_{k})$ for each $\varphi\in C_{\ell, Lip}$. In what follows, such measure families and parametric sets will be fixed for each $k$. 
Denote
$$
\gamma^\theta_{n,k}(n) := \mathsf E^\theta \tilde\chi_n(|X_k|)
\le \mathsf E^\theta 1(|X_k|>n) =: \hat \gamma^\theta_k(n), \quad  \theta\in \Theta_{k},
$$
and for $y>0$,
$$
\gamma^\theta_{n,k}(ny) := \mathsf E^\theta \tilde\chi_n(|X_k|/y)
\le \mathsf E^\theta 1(|X_k|>ny) =: \hat \gamma^\theta_k(ny), \quad  \theta\in \Theta_{k},
$$
or, with $t=ny$ ($y=t/n$),
$$
\gamma^\theta_{n,k}(t) := \mathsf E^\theta \tilde\chi_n(n|X_k|/t)
\le \mathsf E^\theta 1(|X_k|>t) =: \hat \gamma^\theta_{k}(t), \quad  \theta\in \Theta_{k},
$$
and 
$$
\gamma_{n,k}(t) := \mathbb E \tilde\chi_n(n|X_k|/t) =\sup_{\theta\in \Theta_{k}}\mathsf E^\theta \tilde\chi_n(n|X_k|/t)\le \overbrace{ \mathbb E 1(|X_k|>t)}^{=V(|X_k|>t)} =: \hat \gamma_{k}(t), 
$$
and, finally, 
$$
\psi_n(y) := \sum_{k=1}^n y \gamma_k(ny) \equiv \sum_{k=1}^n y\, \mathbb E 1(|X_k|>ny).
$$
NB: Notice that 
$$
\hat \gamma_k((n+1)y) \le \hat \gamma_k(ny).
$$

\begin{theorem}[Nonlinear Kolmogorov's WLLN]\label{thm1}
Let r.v's $(Y_{n,k}, 1\le k\le n)$ be nonlinearly pairwise independent according to the definition \ref{def-indep} for each $n$. Suppose that 
\begin{equation}\label{psin20}
\psi_n(y) \to 0, \; \forall \; y\in [0,1]\quad n\to\infty, 
\end{equation}
and assume that the family of functions $(\psi_n(y), \, y\in [0,1], \, n\ge 1)$ be uniformly integrable.
Then for any $\varepsilon >0$
\begin{equation}\label{Ve}
V(\frac{S_n}{n} \ge \bar \mu_n + \varepsilon) \to 0, \quad n\to\infty, 
\end{equation}
and 
\begin{equation}\label{Ve2}
V(\frac{S_n}{n} \le \underline \mu_n - \varepsilon) \to 0, \quad n\to\infty. 
\end{equation}
Respectively, for any $\varepsilon>0$, 
\begin{equation}\label{ve}
v(\underline \mu_n - \varepsilon < \frac{S_n}{n} < \bar \mu_n + \varepsilon) \to 1, \quad n\to\infty, 
\end{equation}

\end{theorem}
In words, for large values of $n$, the ratio $S_n/n$ will fluctuate between the two moving frontieres $\underline \mu_n$ and $\bar \mu_n$  with a sublinear probability minorated by some sequence of real values approaching one.

\begin{Remark}
Notice that if $V(A) = {\mathbb E} 1(A)\approx 1$, then
$$
v(\Omega \setminus A) = {\cal E}1(\Omega \setminus A) 
= -{\mathbb E} (- 1(\Omega \setminus A)) = 1 - {\mathbb E} (1- 1(\Omega \setminus A)) = 1 - {\mathbb E} 1(A) \approx 0.
$$
Vice versa, if $V(A)\approx 0$, then
$$
v(\Omega \setminus A) = 1 - {\mathbb E} 1(A)   \approx 1.
$$\end{Remark}

\begin{Remark}
A bit more relaxed condition could be used, 
$$
\sum^n\sup_{\theta \in \Theta_{k}} \int_0^1 y  \gamma_k^{\theta}(ny)dy \to 0, 
$$
or, equivalently,??? 
\begin{equation}\label{Psi}
\sum^n \mathbb E\int_0^1 y 1(|X_k| > ny)dy \to 0, \quad n\to\infty.
\end{equation}
However, it is less evident how to verify this condition.
\end{Remark}

\begin{Remark}\label{rem3}
Some extensions are possible for the condition of pairwise independence. 

\begin{itemize}
\item
Instead of the nonlinear pairwise independence of all $(X_k)$, it suffices that for each $n$ the r.v. $(Y_{n,k})$ were pairwise uncorrelated, 
$$
\mathbb E(Y_{n,k} - \tilde \mu_{n,k}) (Y_{n,i} - \tilde \mu_{n,i}) = 0.
$$

\item
Let 
\begin{align*}
&\kappa_{n,i,k}: = \mathbb E(Y_{n,k} - \tilde \mu_{n,k}) (Y_{n,i} - \tilde \mu_{n,i}) 
 \\\\
&= 
\mathbb E(Y_{n,k} - \tilde \mu_{n,k}) (Y_{n,i} - \tilde \mu_{n,i}) - \overbrace{\mathbb E(Y_{n,k} - \tilde \mu_{n,k}) \mathbb E(Y_{n,i} - \tilde \mu_{n,i})}^{=0}.
\end{align*}
What we need instead of the nonlinear pairwise independence is the following assumption, which may also be called  Ces\'aro negative correlation condition, 
\begin{equation}\label{weaknldep}
\frac1{n^2}\left(\sum_{1\le i,k \le n: i\neq k }\kappa_{n,i,k}\right)^+ = o(1), \quad n\to\infty.
\end{equation}
Of course, (\ref{weaknldep}) is  met if $\sum_{1\le i,k \le n: i\neq k }\kappa_{i,k,n} \le 0$ for each $n$, see step 4 in the proof in what follows. 
Here $a^+ = a \,1(a\ge 0)$. 
\end{itemize}

\end{Remark}

\begin{Remark}
For the ``usual'' linear setting in the i.i.d. case and with just one probability measure $\mathsf P$, the combination of conditions on uniform integrability of the family of functions $(\psi_n)$ along with the assumption (\ref{psin20}) is equivalent to the Kolmogorov condition (\ref{ANKcondition}). For the non-i.i.d. case the analogue of this condition was used in  \cite[Theorem 3]{Akhmiarova-V}, where  the r.v's may also be weakly dependent in a certain sense different from (\ref{weaknldep}).
\end{Remark}

\fi


\ifhideen
\section{Proof}
The suggested scheme of the proof is similar to that of the proof of Kolmogorov's classical result (e.g., from \cite[Vol. 2, ch. VII \& ch. XVII]{Feller}). 
At the same time, the nonlinearity clearly brings certain new nuances.

\medskip

\noindent
{\bf I.} The proof of (\ref{Ve}).

\medskip

{\bf 1.} We have, 
\begin{align}\label{VSSprime}
&V(\frac{S_n}{n} \ge \bar \mu_n + \varepsilon) 
= \mathbb E 1(\frac{S_n}{n} \ge \bar \mu_n + \varepsilon)
\,(= \sup_\theta \mathsf P^\theta (\frac{S_n}{n} \ge \bar \mu_n + \varepsilon))
 \nonumber \\\nonumber \\
&
\le  \mathbb E 1(\frac{S'_n}{n} \ge \bar \mu_n + \varepsilon) +  \mathbb E 1(S_n \neq S'_n)
= V(\frac{S'_n}{n} \ge \bar \mu_n + \varepsilon) + V(S_n \neq S'_n).
\end{align}

{\bf 2.} Further, 
\begin{align}\label{SSprime}
&V(S_n \neq S'_n) 
= \mathbb E 1(S_n \neq S'_n)  \le \sum_{k=1}^n \mathbb E 1(Y_k \neq X_k)
 \nonumber \\\nonumber \\
&  
= \sum_{k=1}^n  V(|X_k| >n) = \sum_{k=1}^n \hat \gamma_k(n) = \psi_n(1) \to 0, \quad n\to\infty, 
\end{align}
by the assumption (\ref{psin20}).

{\bf 3.} Further, 
\begin{align}\label{VSprime}
&V(\frac{S'_n}{n} \ge \bar \mu_n + \varepsilon) = V(S'_n\ge n\bar \mu_n + n\varepsilon) = V(S'_n\ge  n\varepsilon+ \sum_{k=1}^n \tilde \mu^+_{n,k})
 \nonumber \\ \nonumber \\
& \le V(|S'_n - \sum_{k=1}^n \tilde \mu^+_{n,k}|\ge  n\varepsilon) \le \frac1{n^2\varepsilon^2} \mathbb E (S'_n -  \sum_{k=1}^n \tilde \mu^+_{n,k})^2, 
\end{align}
by virtue of Markov's inequality (see the lemma \ref{Markov}). 

{\bf 4.}
Let us estimate the value $ \mathbb E (S'_n -  \sum_{k=1}^n \tilde \mu^+_{n,k})^2$. We have, 
\begin{align*}
& \mathbb E (S'_n -  \sum_{k=1}^n \tilde \mu^+_{n,k})^2 = \mathbb E (\sum_{k=1}^n (Y_{n,k} - \tilde \mu^+_{n,k}))^2 
 \\\\
& = \mathbb E \left(\sum_{k=1}^n (Y_{n,k} - \tilde \mu^+_{n,k})^2 + \sum_{i,k: i\neq k} (Y_{n,k} - \tilde \mu^+_{n,k}) (Y_{n,i} - \tilde \mu^+_{n,i})\right)
\\\\
& \le  \sum_{k=1}^n \mathbb E (Y_{n,k} - \tilde \mu^+_{n,k})^2 + \sum_{i,k: i\neq k}  \mathbb E(Y_{n,k} - \tilde \mu^+_{n,k}) (Y_{n,i} - \tilde \mu^+_{n,k}).
\end{align*}
By the accepted definition of pairwise nonlinear r.v. independence, and because of the definition of the values $\tilde \mu^+_{n,k}$, for $i\neq k$, we have, 
\begin{equation}\label{nlindep}
\mathbb E(Y_{n,k} - \tilde \mu^+_{n,k}) (Y_{n,i} - \tilde \mu^+_{n,i}) = \mathbb E(Y_{n,k} - \tilde \mu^+_{n,k}) \mathbb E(Y_{n,i}- \tilde \mu^+_{n,i}) = 0.
\end{equation}
{\em NB: This is where the condition (\ref{weaknldep}) from the remark \ref{rem3} may be applied instead of the equality (\ref{nlindep}).  It would suffice for the rest of the proof.}\\
In any case, it follows from (\ref{nlindep}) that
\begin{align*}
& \mathbb E (S'_n -  \sum_{k=1}^n \tilde \mu^+_{n,k})^2 \le  \sum_{k=1}^n \mathbb E (Y_{n,k} - \tilde \mu^+_{n,k})^2 =  \sum_{k=1}^n \mathbb E (Y_{n,k}^2 - 2\tilde \mu^+_{n,k} Y_{n,k} + (\tilde \mu^+_{n,k})^2)
 \\\\
& \le \sum_{k=1}^n (\mathbb E Y_{n,k}^2 + \mathbb E (-2\tilde \mu^+_{n,k} Y_{n,k}) + \mathbb E (\tilde\mu^+_{n,k})^2) 
 \\\\
& = \sum_{k=1}^n \mathbb E Y_{n,k}^2 + \sum_{k=1}^n \mathbb E (-2\tilde \mu^+_{n,k} Y_{n,k}) +  \sum_{k=1}^n  (\tilde\mu^+_{n,k})^2.
\end{align*}

{\bf 5.} Let us estimate each sum in the last expression. We have, 
$$
\mathbb E Y_{n,k}^2 = \sup_{\theta\in \Theta_{k}} \mathsf P^\theta Y_{n,k}^2.
$$
Denote 
$$
F^\theta_k(x) := \mathsf P^\theta (X_k\le x), \quad \theta\in \Theta_{k}.
$$
Since each $\mathsf P^\theta$ is a sigma-additive probability measure, the function  $F^\theta_k$ for each $\theta\in \Theta_{k}$ is a usual distribution function of the r.v. $X_k$ with respect to this probability.
Following  \cite[Vol. 2, Chapter VII, \S 7, formula (7.7)]{Feller}, let us introduce the function 
$$
\sigma^\theta_k(t):= \frac1{t} \int_{-t}^t x^2 dF^\theta_{X_k}(x), \quad \theta\in \Theta_{k},
$$
and transform it by integrating by parts, as suggested in \cite{Feller}, 
$$
\sigma^\theta_k(t)= \frac1{t} \int_{-t}^t x^2 dF^\theta_{X_k}(x) = 
- t \hat \gamma^\theta_k(t) + \frac2{t} \int_0^t x \hat \gamma^\theta_k(x)dx \le \frac2{t} \int_0^t x \hat \gamma^\theta_k(x)dx.
$$
Then we estimate, using change of variables, 
\begin{align}\label{I}
&\frac1{n^2} \sum_{k=1}^n \mathbb E (Y_{n,k}^2) = \frac1{n^2} \sum_{k=1}^n \sup_{\theta \in \Theta_{k}} \mathsf E^\theta (Y_{n,k}^2) 
= \frac1{n^2} \sum_{k=1}^n \sup_{\theta \in \Theta_{k}} \mathsf E^\theta (X_{k}^2 \chi^2_{n}(|X_{k}|)) 
 \nonumber \\ \nonumber \\
&\! \le \frac1{n^2} \sum_{k=1}^n \sup_{\theta \in \Theta_{k}} \mathsf E^\theta (X_{k}^2 1(|X_{k}| \le n+1))
=  \frac1{n^2} \sum_{k=1}^n \sup_{\theta \in \Theta_{k}} \int\limits_{-n-1}^{n+1} x^2 dF^\theta_{X_k}(x)
 \nonumber \\ \nonumber \\
&\! =   \!\frac{n+1}{n^2} \sum_{k=1}^n \sup_{\theta \in \Theta_{k}} \sigma^\theta_k(n+1)
\!\le \!  \frac{2}{n^2} \sum_{k=1}^n \sup_{\theta \in \Theta_{k}}  \int\limits_0^{n+1} x \hat \gamma^\theta_k(x)dx 
 \nonumber \\\nonumber \\
&
\!\le\! \frac2{n^2} \sum_{k=1}^n  \int\limits_0^{n+1} x \underbrace{\sup_{\theta \in \Theta_{k}}\hat \gamma^\theta_k(x)}_{=\hat \gamma_k(x)}dx
\stackrel{y=x/(n+1)}= 2 \frac{(n+1)^2}{n^2}\int\limits_0^1  \sum_{k=1}^n y \hat \gamma_k((n+1)y) dy 
 \nonumber \\\nonumber \\
&
\le 2 \frac{(n+1)^2}{n^2}\int\limits_0^1  \sum_{k=1}^n y \hat \gamma_k(ny) dy 
= 2 \left(1+\frac1{n}\right)^2 \int_0^1 \psi_n(y)dy \to 0, \quad n\to\infty, 
\end{align}
the latter convergence 
due to the assumptions (\ref{psin20}) and UI of the family $(\psi_n(y), \, y\in [0,1])$.


\medskip

{\bf 6.} Let us evaluate the term $\sum_{k=1}^n  (\tilde\mu^+_{n,k})^2$. We have, 
\begin{align}\label{II}
& \frac1{n^2} \sum_{k=1}^n  (\tilde\mu^+_{n,k})^2 = \frac1{n^2} \sum_{k=1}^n  (\mathbb E (Y_{n,k}))^2 = \frac1{n^2} \sum_{k=1}^n  (\sup_{\theta \in \Theta_{k}} \mathsf E^\theta (Y_{n,k}))^2
 \nonumber \\ \nonumber \\
& =  \frac1{n^2} \sum_{k=1}^n  (\sup_{\theta \in \Theta_{k}} \mathsf E^\theta (X_k 1(|X_k|\le n)))^2 \le \frac1{n^2} \sum_{k=1}^n  (\sup_{\theta \in \Theta_{k}} \mathsf E^\theta (|X_k| 1(|X_k|\le n)))^2 
 \nonumber \\ \nonumber \\
& \stackrel{\text{CBS}} \le \frac1{n^2} \sum_{k=1}^n  \sup_{\theta \in \Theta_{k}} \mathsf E^\theta (|X_k|^2 1(|X_k|\le n)) = \frac1{n^2} \sum_{k=1}^n \mathbb E (Y_{n,k}^2) \to 0, \; n\to\infty, 
\end{align}
due to what was just established in (\ref{I}). Here ``CBS'' stands for Cauchy -- Buniakowskii -- Schwarz' inequality applied for each fixed $\theta$ and probability~$\mathsf P^\theta$.

\medskip

{\bf 7.} Finally, to estmate the remaining sum $\sum_{k=1}^n \mathbb E (-2\tilde \mu^+_{n,k} Y_{n,k})$, we write, 
\begin{align}\label{III}
& \frac1{n^2} \sum_{k=1}^n \mathbb E (-2\tilde \mu^+_{n,k} Y_{n,k}) 
\le  \frac2{n^2} \sum_{k=1}^n |\tilde \mu^+_{n,k}| \mathbb E |Y_{n,k}| 
=  \frac2{n^2} \sum_{k=1}^n (|\mathbb E Y_{n,k}|) (\mathbb E |Y_{n,k}|) 
 \nonumber\\\nonumber\\
& \le \frac2{n^2} \sum_{k=1}^n (\mathbb E |Y_{n,k}|)^2 =  \frac2{n^2} \sum_{k=1}^n  (\sup_{\theta \in \Theta_{k}}\mathsf E^\theta 1(|X_{k}| 1(|X_{k}|\le n)))^2 
 \nonumber\\\nonumber\\
& \stackrel{\text{CBS}}\le \frac2{n^2} \sum_{k=1}^n  (\sup_{\theta \in \Theta_{k}}\mathsf E^\theta 1(|X_{k}|^2 1(|X_{k}|\le n))) 
=  \frac2{n^2} \sum_{k=1}^n \mathbb E (Y_{n,k}^2) \to 0, \; n\to\infty, 
\end{align}
again by virtue of (\ref{I}). 

\medskip

Overall, combining all bounds (\ref{I}) -- (\ref{III}), we get the claim (\ref{Ve}) from (\ref{VSSprime}), (\ref{SSprime}), and (\ref{VSprime}), 
as required. 

\medskip


\noindent
{\bf II.} The proof of (\ref{ve}) is similar. Indeed, we have, 

\begin{align}\label{vSSprime}
&V(\frac{S_n}{n} \le \underline \mu_n - \varepsilon) = \mathbb E 1(\frac{-S_n}{n} \ge - \underline \mu_n + \varepsilon)) 
 \nonumber \\\nonumber \\
&\le   \mathbb E 1(\frac{-S'_n}{n} \ge -\underline \mu_n + \varepsilon) +  \mathbb E 1(-S_n \neq -S'_n)
 \nonumber \\\nonumber \\
&= V(\frac{-S'_n}{n} \ge -\underline \mu_n + \varepsilon) + V(S_n \neq S'_n).
\end{align}
The rest follows from the previous calculus applied to the sequences $(-X_k)$ and the ``lower nonlinear mean values'' $\tilde \mu^-_{n,k}$.

\medskip

\noindent
{\bf III.} Let $\displaystyle A(n,\varepsilon) := (\omega: \frac{S_n}{n} \ge \bar\mu_n + \varepsilon, \;\text{or} \;\frac{S_n}{n} \le \underline \mu_n - \varepsilon)$. By the sublinear expectation rules, we have, 
$$
V(A(n,\varepsilon)) \le V(\frac{S_n}{n} \ge \bar\mu_n + \varepsilon) + V(\frac{S_n}{n} \le \underline \mu_n - \varepsilon) \to 0, \;\; n\to\infty.
$$
Hence, by the definition of the operator $v$, we obtain, 
$$
v(\underline \mu_n - \varepsilon < \frac{S_n}{n} <  \bar\mu_n + \varepsilon) = v(\Omega \setminus A(n,\varepsilon)) = 1 - V(A(n,\varepsilon)) \to 1, \quad n\to\infty.
$$
The theorem is proved. \hfill $\square$

\medskip

\begin{Remark}
Notice that although we used the representation of the sublinear expectation $\mathbb E$ for each $X_k$ as $\mathbb E g(X_k) = \sup_{\theta\in \Theta_{k}} \mathsf  E^\theta\ g(X_k)$, eventually all assumptions were in terms of $\mathbb E$, not in terms of individual expectations $\gamma^\theta_{n,k}$, et al. So, any particular choice of the family $\mathsf P^\theta, \theta\in \Theta_{k}$ is of no importance.

\end{Remark}

\fi


\ifhideen
\section*{Acknowledgements}
\noindent
For both authors this study was funded by the Theoretical Physics and Mathematics Advancement Foundation ``BASIS''. 
\fi

\ifhideen
\section*{Co-authors' contributions}
\noindent
The first author made the first version of the paper and the calculus under the assumption in the example \ref{Ex1}. 
The second author suggested the hypothesis, performed a general supervision and made a final polishing under the assumptions in terms of the abstract nonlinear expectation $\mathbb E$. 
\fi

\ifhideen

\fi

\ifhideru

\fi

\end{document}